\documentclass[12pt]{article}
\usepackage{amssymb}
\usepackage{amsmath}
\usepackage{graphicx}
\usepackage{acronym}
\usepackage{amsfonts}
\usepackage{color}
\usepackage{colordvi}
\newcommand{\qed}{\bull \medskip}

\newtheorem{theorem}{Theorem} 
\newtheorem{corollary}{Corollary}
\newtheorem{proposition}{Proposition}
\newtheorem{lemma}{Lemma}
\newtheorem{definition}{Definition}
\def\be{\begin{Example}}
\def\ee{\end{Example}}
\def\bt{\begin{theorem}}
\def\et{\end{theorem}\bigskip}
\def\bl{\begin{Lemma}}
\def\el{\end{Lemma}\bigskip}
\def\ep{\end{Proposition}\bigskip}
\def\bp{\begin{Proposition}}
\def\bd{\begin{definition}}
\def\ed{\end{definition}}

\def\qed{\hfill {$ \Box $} \medskip}
\def\proof{\noindent\bf Proof. \hspace{4mm}\rm}

\newcommand{\A}{{\cal A}}

\newcommand{\C}{{\Bbb C}}
\newcommand{\Z}{{\Bbb Z}}
\newcommand{\R}{{\Bbb R}}
\newcommand{\PP}{{\Bbb P}}
\newcommand{\x}{{\bf x}}

\begin{document}
\title{\bf E-Characteristic Polynomials of Tensors\thanks
 {To appear in: Communications in Mathematical Sciences.  This work is supported by the Research Grant Council of Hong Kong, Project Numbers: PolyU 501808, 501909
 and 502510.}}
\author{An-Min Li \footnote{\small School of Mathematics, Sichuan University, China.
E-mail: math\_li@yahoo.com.cn}
 \ \
Liqun Qi \footnote{\small Department of
Applied Mathematics, The Hong Kong Polytechnic University, Hung Hom,
Kowloon, Hong Kong. E-mail: maqilq@polyu.edu.hk}
\ \
Bin Zhang \footnote{\small School of Mathematics, Sichuan University, China.
E-mail: zhangbin@scu.edu.cn}}
\date{\today}
\maketitle

{\large
\begin{abstract} In this paper, we show that the coefficients of the
E-characteristic polynomial of a tensor are orthonormal invariants
of that tensor.  When the dimension is $2$, some simplified formulas
of the E-characteristic polynomial are presented.    A resultant
formula for the constant term of the E-characteristic polynomial is
given. We then study the set of tensors with infinitely  many
eigenpairs  and the set of irregular tensors, and prove both the
sets have codimension 2 as subvarieties in the projective space of
tensors. This makes our perturbation method workable. By using the
perturbation method and exploring the difference between
E-eigenvalues and eigenpair equivalence classes, we present a simple
formula for the coefficient of the leading term of the
E-characteristic polynomial, when the dimension is $2$.
\\
{\bf Keywords:}  E-Eigenvalues, tensors, E-characteristic
polynomials, eigenpair equivalence class, irregularity

\end{abstract}

\section{Introduction} Eigenvalues of higher order tensors were introduced in
2005 \cite{qi, li} and have attracted much attention in the
literature and found applications in science and engineering.   The
E-eigenvalues of a tensor were introduced in \cite{qi, qi07}.
E-eigenvalues are invariant under orthonormal coordinate changes \cite{qi07}.   They have applications in determining
positive definiteness of a multivariate system \cite{qi, qww}, best
rank-one approximation \cite{qi, qi11, qww, zqy}, magnetic resonance
imaging \cite{bv, qyw}, spectral hypergraph theory \cite{hq} and
symmetric embedding \cite{rv}, and possess links with geometry
\cite{ba, ba09, ba11, ba10, qi06}. The numbers of E-eigenvalues and
E-eigenvectors have been studied in \cite{cpz, nqww}.   For a real
tensor, an E-eigenvalue with a real E-eigenvector is called a
Z-eigenvalue. Z-eigenvalues are real and play the main role in the
above applications. Algorithms for finding Z-eigenvalues were
studied in \cite{km, qww}.

An $n$-dimensional $m$th order hypermatrix  $A$ with entries in $\C $ is a map
from $\{1,\cdots, n\}^m \to \C$. We use notation $A= ( a_{i_1i_2\cdots i_m})$ to denote such a hypermatrix.

In physics \cite{hs, wg}, a tensor is a physical quantity.    In geometry, a tensor can be regarded as a multi-linear function in a Hilbert
space, see discussion in Section 3.   In such applications, in an orthonormal coordinate system,
a tensor is expressed by a hypermatrix.   It is important to determine parameters and properties of
such a tensor, which are invariant under orthonormal coordinate changes.    In Theorem 1 of Section of this paper,
we will show that the coefficients of the E-characteristic polynomial of a tensor are orthonormal invariants of that
tensor.   For that discussion, we need to distinguish a tensor $\A$ and its hypermatrix $A$ in a coordinate system.

For some other applications, such as determining
positive definiteness of a multivariate system, best
rank-one approximation,  spectral hypergraph theory and
symmetric embedding, there are no coordinate systems involved.
People just regard a hypermatrix $A$ as a tensor.   Hence, in the
other parts of this paper, a tensor simply means a hypermatrix.

For an $m$th order tensor $A$ and a vector $x \in \C^n$, we define
$A x^{m-1}$ to be the vector in $\C^n$ with its $i$th component as
$$(A x^{m-1})_i =\sum\limits_{i_2,\cdots,i_m=1}^n
a_{ii_2\cdots i_m}x_{i_2}\cdots x_{i_m}.$$   If $\lambda\in \C$ and
$x \in \C^n$ satisfy
\begin{equation} \label{eig}
\left\{{A x^{m-1} = \lambda x, \atop x^Tx = 1,}\right.
\end{equation}
then $\lambda$ is an E-eigenvalue of $A$ and $x$ is an E-eigenvector
of $A$.   If $A$ and $x$ are real, then $\lambda$ is also real
\cite{qi}.  In this case, $\lambda$ is called a Z-eigenvalue of $A$
and $x$ is called a Z-eigenvector of $A$.

E-eigenvalues and Z-eigenvalues are invariant under orthonormal
coordinate changes \cite{qi}. An even-order real symmetric
tensor always has Z-eigenvalues \cite{qi}.   It is positive definite
(semi-definite) if and only if all of its Z-eigenvalues are positive
(nonnegative) \cite{qi}. The Z-eigenvalue with the largest absolute
value and its corresponding Z-eigenvector form the best rank-one
approximation of a real symmetric tensor \cite{qi, qww}.   The
modifier ``E-'' stands for Euclidean, as (\ref{eig}) implies that
the Euclidean norm of $x$ is $1$ if $x$ is real.   The modifier
``Z-'' names after Prof. Shuzi Zhou, who suggested (\ref{eig}) to
the author of \cite{qi}.  Prof. Zhou departed by cancer in 2009. The
name ``Z-eigenvalue'' has already been used widely in the literature
\cite{ba, ba09, ba11, ba10, bv, cpz, hq, km, nqww, qi, qi06, qi07,
qi11, qww, qyw, za, zqy}. Hence, we keep to use this name.

E-characteristic polynomials were introduced in \cite{qi, qi07}, and
discussed in \cite{nqww, cs}.   Until now, there are no other papers
containing discussion on E-characteristic polynomials. We now survey
the definitions and properties of E-characteristic polynomials in \cite {cs, nqww, qi, qi07}

\begin{definition} \label{d1}
Let $A = ( a_{i_1i_2\cdots i_m})$ be an $n$-dimensional $m$th order
tensor, if there exists $x \in \C^n \setminus \{ 0 \}$ such that
\begin{equation} \label{eig2}
\left\{{A x^{m-1} = 0, \atop x^Tx = 0,}\right.
\end{equation}
then $A$ is called irregular.  Otherwise, $A$ is called regular
\cite{qi, qi07}.

When $m$ is even, the E-characteristic polynomial $\psi_A(\lambda)$
is defined \cite{qi} as
\begin{equation} \label{ch1}
\psi_A(\lambda) = {\rm Res}_x\left(A x^{m-1} - \lambda \left(x^Tx\right)^{{m-2 \over 2}}x\right) = \sum_{j=0}^{h(m, n)} a_j(m,
n)\lambda^j,
\end{equation}
where the second equality is the expansion of the resultant of (\ref {eig2}) in terms of $\lambda$, $h(m,n)$ is the highest power with generically non-zero coefficient and $a_j(m,n)$'s are polynomials in entries of $A$.

When $m$ is odd, the E-characteristic polynomial $\psi_A(\lambda)$
is defined \cite{qi07} as
\begin{equation} \label{ch2}
\psi_A(\lambda) = {\rm Res}_{(x, x_0)}\left({A x^{m-1} - \lambda
x_0^{m-2}x \atop x^Tx - x_0^2}\right) = \sum_{j=0}^{h(m, n)}
a_j(m, n)\lambda^{2j},
\end{equation}
where $2h(m,n)$ is the highest power with generically non-zero coefficient.
\end{definition}

\begin {lemma} For a regular tensor, E-characteristic polynomial can be defined as  $$Res _{(x,x_0)}\left (\begin {array}{c}Ax^{m-1}-\lambda x_0^{m-2}x\\x^Tx-x_0^2\end{array}\right ),
$$
for all $m$.
\end{lemma}
\proof For a regular tensor and even $m$, the systems
$$\left \{\begin {array}{c} Ax^{m-1}-\lambda (x_0)^{m-2}x=0\\x^Tx-x_0^2=0\end{array}\right .
$$
and
$$Ax^{m-1}-\lambda (x^Tx)^{\frac {m-2}2}x=0
$$
are equivalent. \qed

According to Theorem 4 of \cite{qi07}, any E-eigenvalue $\lambda$ of
$A$ must be a root of $\psi_A(\lambda)$; if $A$ is regular, then a
complex number $\lambda$ is an E-eigenvalue of $A$ if and only if it
is a root of $\psi_A(\lambda)$.   It was proved in \cite{nqww} that
when $m$ is even,
\begin{equation} \label{equ5}
h(m, n) = {(m-1)^n -1 \over m-2} = \sum_{i=0}^{n-1} (m-1)^i.
\end{equation}

The definition (\ref{eig}) is not a strict extension of the
classical definition for eigenvalues of a square matrix.  The main
point is that it excludes complex eigenvalues whose eigenvectors
satisfy $x^Tx = 0$ but $x \not = 0$ \cite{qi07}.   Recently,
based upon this observation, Cartwright and Sturmfels \cite{cs}
introduced equivalence classes of eigenpairs. If $\lambda\in \C$ and
$x \in \C^n \setminus \{ 0 \}$ satisfy
\begin{equation} \label{eig1}
A x^{m-1} = \lambda x,
\end{equation}
then $(\lambda, x)$ is called an eigenpair of $A$.  Two eigenpairs
$(\lambda, x)$ and $(\lambda^\prime, x^\prime)$ are considered to be
equivalent if there is a complex number $t \not = 0$ such that
$\lambda^\prime = t^{m-2}\lambda$ and $x^\prime = tx$.
Cartwright and Sturmfels proved that (\ref {equ5}) is true for all
$m$ when counting the number of generalized
eigenpairs without any normalization restriction \cite{cs}.

An important property of a second-order tensor is that not only its
eigenvalues but also the coefficients of its characteristic
polynomial are invariants of that tensor. It was proved in \cite{qi,
qi07} that E-eigenvalues of a higher order tensor are invariants of
that tensor.   This implies that the coefficients of the
E-characteristic polynomial, divided by the first nonzero
coefficient, are invariants of that tensor, see Theorem 5 of
\cite{qi07}.   There was no existing result about the invariance of
the coefficients of the E-characteristic polynomial.

In \cite{cs}, it was given that
$$a_0(3, 2) = \left({\rm Res}_x\left(A x^2\right)\right)^2,$$
and mentioned in general that $a_0(m, n)$ is a power of ${\rm
Res}_x\left(A x^{m-1}\right)$, without specifying the value of the
power or giving a proof.      When $m=3$ and $n=2$, they also
present the coefficient of the leading term as the negative of the
sum of two squares of linear combinations of the entries of $A$.

In this paper, we explore more in this direction.  The structure of
the paper is as following.

Since resultants are main tools to study the E-characteristic polynomials, we recall some facts about resultants in Section 2. In Section 3, we prove that all the coefficients of the E-characteristic polynomial
are orthonormal invariants of a tensor. We then show that when $m$ is even, $a_0(m, n)$ is Res$_x(A x^{m-1})$,
and when $m$ is odd, it is the square of Res$_x(A x^{m-1})$.

In Section 4, we give two simplified determinantal formulas
for the E-characteristic polynomial when $n=2$.    In the even case,
this is a $(2m-2) \times (2m-2)$ determinant.   In the odd case,
this is a $(3m-4) \times (3m-4)$ determinant.   The parameter
$\lambda$ only appears in the first $m$ rows of the determinants.

We then address the leading coefficients of E-characteristic
polynomials. Two sets of tensors are easy to handle: i) the set of
tensors with infinitely many eigenpairs, ii) the set of regular
tensors.   Actually, the E-characteristic polynomials of tensors in
the first set are identically zero; while a complex number is an
E-eigenvalue of a regular tensor if and only if it is a root of the
E-characteristic polynomial of that tensor. For tensors with
finitely many eigenpairs, we apply the perturbation method, to
approximate them by regular tensors with finitely many eigenpairs.
Therefore we have to make sure that the perturbation is possible.
Sections 5 and 6 serve this purpose. In Section 5, we prove that the
set of tensors with infinitely many eigenpairs has codimension at
least 2 in the projective space of tensors. This shows it is
possible to approximate a tensor in a hypersurface in the projective
space of tensors by tensors with finitely many eigenpairs. We
show in Section 6 that the set of irregular tensors has codimension
at least 2 in the projective space of tensors too. Therefore, the
union of the set of tensors with infinitely many eigenpairs and the
set of irregular tensors has codimension at least 2 in the
projective space of tensors, which make the approximation to a
tensor on a hypersurface by by regular tensors with finitely many
eigenpairs possible.

Finally, in Section 7, we introduce
the deficit system to explore the difference  of E-eigenvalues and
eigenpair equivalence classes of Cartwright and Sturmfels \cite{cs}.
Based on this analysis, we show that when $n=2$, the leading
coefficient of the E-characteristic polynomial is the $({m-2 \over
2})$th power of the sum of two squares when $m$ is even, and the
negative of the $(m-2)$th power of such a sum when $m$ is odd.

\section {Resultants}
\label {s:resultant}

First we recall some results for resultants since our discussion heavily depends on the properties of resultants, according to \cite{clo} and \cite{gkz}.

For fixed positive integers $d_1, \cdots, d_n$, and each pair $(i,\alpha )$, where $i=1, \cdots n$, $\alpha =(\alpha _1, \cdots , \alpha _n)\in \Z _{\ge 0}^n$ with $|\alpha |=\sum \alpha _j =d_i$, we introduce a variable $u_{i,\alpha}$. We have the following results about resultants.

\begin {proposition} There is a unique polynomial $Res_x=Res_{d_1, \cdots , d_n}\in \Z [u_{i,\alpha}]$ with the following properties:
\begin {itemize}
\item [a.] If $F_1, \cdots , F_n \in \C [x_1, \cdots, x_n]$ are homogeneous of degrees $d_1, \cdots, d_n$, then the equations $F_1=\cdots =F_n=0$ have a solution in $\C \PP ^{n-1}$ if and only if $Res(F_1, \cdots , F_n)=0$, where the notation $Res (F_1, \cdots, F_n)$ means $Res (c_{i,\alpha})$, if $$F_i=\sum_{|\alpha |=d_i} c_{i,\alpha }x^\alpha.
$$
\item [b.] $Res(x_1^{d_1}, \cdots, x_n^{d_n})=1$.
\item [c.] $Res$ is irreducible in $\C [u_{i,\alpha}]$.
\end{itemize}
\end{proposition}

\begin{proposition} \label{t01} Res$(F_1, \cdots, F_n)$ is a homogeneous polynomial in
the coefficients of $F_i$, with degree $d_1\cdots
d_{i-1}d_{i+1}\cdots d_n$, for $i = 1, \cdots, n$.
\end{proposition}

\begin{proposition} \label{res02} If $F_j$'s are homogeneous of degree $d$, and $G_i=\sum a_{ij} F_j$, where $(a_{ij})$ is an invertible matrix, then
$$Res (G_1, \cdots, G_n)=det(a_{ij})^{d^{n-1}}Res(F_1, \cdots, F_n).
$$
\end{proposition}

\begin{proposition} \label{res03} If $F_j=F'_jF^{\prime \prime }_j$ is a product of homogeneous polynomials,
$$Res (F_1, \cdots, F_j, \cdots, F_n)=Res (F_1, \cdots, F'_j, \cdots, F_n)Res (F_1, \cdots, F^{\prime\prime} _j, \cdots, F_n).
$$
\end{proposition}

For an $n$-dimensional vector space $V$ (over $\R $ or $\C$), let $L: V\to V$ be a linear transformation, it induces a natural map $L^*$ on functions on $V$. Fix a basis $\{e_1, \cdots, e_n\}$ of $V$, $x_i(\sum c_je_j)=c_i$ defines a linear function $x_i$ on $V$, for $i=1,\cdots n$. If the  transformation matrix of $L$ under this basis is $(a_{ij})$, then $L^*x_i=\sum a_{ij}x_j$. So for a polynomial $F$ in variables $x_1,\cdots x_n$, it is a function on $V$, and $L^*F$ is a polynomial in variables $x_1,\cdots x_n$ obtained from $F$ by substituting $x_i$ with $\sum a_{ij}x_j$.

\begin {corollary}  For an invertible linear transformation $L: V\to V$,
$$Res (L^*F_1, \cdots, L^*F_n)=det(a_{ij})^{d_1\cdots d_n}Res (F_1, \cdots F_n)
$$
\end{corollary}
\proof Let $$F_i=\sum_{|\alpha |=d_i} u_{i,\alpha }x^\alpha.
$$
Then $Res (F_1, \cdots, F_n)$ is an irreducible polynomial in variables $u_{i,\alpha}$, and by definition, if we denote
$$L^*F_i=\sum_{|\alpha |=d_i} v_{i,\alpha }x^\alpha,
$$
then $v_{i,\alpha}$ is a linear combination of $u$-variables $u_{j,\beta}$,
and $Res (L^*F_1, \cdots, L^*F_n)$ is obtained by substitute $u_{i,\alpha}$ by $v_{i,\alpha}$, so it is a polynomial in variables $u_{i,\alpha}$ and of the same degree as $Res (F_1, \cdots, F_n)$.

Since the matrix $(a_{ij})$ is invertible, the system $F_1=\cdots=F_n=0$ has a nontrivial solution if and only if the system $L^*F_1=\cdots=L^*F_n=0$ has a nontrivial solution, therefore $Res (F_1, \cdots, F_n)$ is a divisor of $Res (L^*F_1, \cdots, L^*F_n)$. We already know they have the same degree, so
$$Res(L^*F_1, \cdots, L^*F_n)=c\ Res(F_1, \cdots, F_n),
$$
where $c$ is a constant.

We can determine this constant by checking the case $F_i=x_i^{d_i}$,
$i=1, \cdots, n$. By definition
$$L^*x_i=\sum_j a_{ij}x_j.
$$
Now by Proposition \ref {res03} and the fact that $Res_{1,\cdots,
1}$ is the standard determinant,
$$\begin {array}{rl}Res (L^*F_1, \cdots, L^*F_n)=&Res ((\sum a_{1j}x_j)^{d_1}, \cdots ,(\sum a_{nj}x_j)^{d_n})\\
=&Res (\sum a_{1j}x_j, (\sum a_{2j}x_j)^{d_2}\cdots ,(\sum a_{nj}x_j)^{d_n})\\
&\times Res ((\sum a_{1j}x_j)^{d_1-1}, \cdots ,(\sum a_{nj}x_j)^{d_n})\\
=&\cdots\\
=&Res (\sum a_{1j}x_j, (\sum a_{2j}x_j)^{d_2}\cdots ,(\sum a_{nj}x_j)^{d_n})^{d_1}\\
=&\cdots\\
=&Res (\sum a_{1j}x_j, \sum a_{2j}x_j, \cdots ,\sum a_{nj}x_j)^{d_1\cdots d_n}\\
=&det(a_{ij})^{d_1\cdots d_n}.\end{array}
$$
Hence,  $$c=det(a_{ij})^{d_1\cdots d_n}.
$$
\qed

\section{Orthonormal Invariance of the Coefficients}
\label {s:invariance}

The E-characteristic polynomial $\psi_A(\lambda)$ is as defined in
the introduction.    Let $a_{\bar h(m, n)}(m, n)$ be the first
nonzero coefficient of $\psi_A(\lambda)$.  Then, for $j = 0, \cdots,
\bar h(m, n) - 1$, according to the relations between roots and
coefficients of a one-variable polynomial,
$$s_{\bar h(m, n) - j} \equiv (-1)^{\bar h(m, n) - j}a_j(m, n)/a_{\bar h(m, n)}(m, n)$$
is the sum of all possible products of $\bar h(m, n) - j$
E-eigenvalues, hence an elementary symmetric function of the
E-eigenvalues.

In the following discussion, we distinguish a tensor $\A$ and its hypermatrix $A$ in
a coordinate system.

Let $V$ be an n-dimensional vector space over $\R$ with a Euclidean
inner product $<,>$. Let  $\A$ be an  $m$th order covariant
tensor on $V$, that is,
$$\A: V^m \rightarrow \R$$
is a multi-linear function. In particular, a first order covariant
tensor is called a covariant vector, denoted by $\x$. In practice, what we consider is complexified tensor, i.e., to view them as multiple-linear functions from $V^m \otimes _\R  \C \to \C$.

We choose an orthonormal frame $E = \{ e_1,...,e_n\}$ in $V$, denote
$$\A(e_{i_1},...,e_{i_m}):= a_{i_1i_2\cdots i_m},\;\; \x(e_i)=x_i.$$
Then $( a_{i_1i_2\cdots i_m})$ (resp. $x_i$) is the hypermatrix of
$\A$ (resp. $\x$) with respect to the frame $E$.  If all
$a_{i_1i_2\cdots i_m}$, for $i_1, \cdots, i_m = 1, \cdots, n$, are
real, then we say that $\A$ is real. If all $x_i$, for $i=1, \cdots,
n$, are real, then we say that $\x$ is real.

If we choose another orthonormal frame $\tilde E = \{
\tilde{e}_1,...,\tilde{e}_n\}$, that is
$$\tilde{e}_i = \sum C_i^j e_j$$
where the matrix $(C_i^j) \in O(n)$ is the transformation matrix, then the hypermatrices of $\A$ and $\x$
with respect to $\tilde{e}_1,...,\tilde{e}_n$ are given by
\begin{equation} \label{e1}
\tilde{a}_{i_1i_2\cdots i_m}= \sum
C_{i_1}^{j_1}...C_{i_m}^{j_m}a_{j_1j_2\cdots j_m},
\end{equation}
\begin{equation} \label{e2}
\tilde{x}_i = \sum C_i^jx_j.
\end{equation}

\noindent {\bf Remark. } For a tensor, the hypermatrices w.r.t.
different orthonormal frames are orthogonally similar.

\noindent {\bf Remark. } Since we only consider tensors under orthonormal
transformations, we do not distinguish covariant tensors and
contravariant tensors.

For a complex $m$th order covariant tensor $\A $ and a covariant vector $\x $, we can define a covariant vector $\A \x^{m-1}$ by
\begin{equation} \label{e4}
\A \x^{m-1}(e_i) = \sum_{i_2, \cdots, i_m=1}^n \A(e_i, e_{i_2}, \cdots,
e_{i_m})\x(e_{i_2})\cdots \x(e_{i_m}),
\end{equation}
for $i= 1, \cdots, n$, and any orthonormal
frame $E = \{ e_1,...,e_n\}$ in $V$.  For any covariant vector $\x $, define
\begin{equation} \label{e5}
\x^T \x = \sum_{i=1}^n \left(\x(e_i)\right)^2,
\end{equation}
for any orthonormal frame $E = \{ e_1,...,e_n\}$ in $V$. By
(\ref{e1}) and (\ref{e2}), we see that the definitions (\ref{e4}) and (\ref{e5}) are well-defined, independent from the
frame $E$.

For a complex $m$th order covariant tensor $\A $ , if $\lambda\in \C$ and non-zero covariant vector $\x $ satisfy
$$
\left\{{\A \x^{m-1} = \lambda \x, \atop \x^T \x = 1,}\right.
$$
then $\lambda$ is an E-eigenvalue of $\A$ and $\x$ is an
E-eigenvector of $\A$. By the above discussion, we see that
$\lambda$ is invariant, i.e., it is independent of the choice of the
frames.

We say that a fact or quantity is {\bf orthonormal invariant}, if it is invariant under the changes of orthonormal frames.

Fix an orthonormal frame $E = \{ e_1,...,e_n\}$ in $V$, let the hypermatrix $A=( a_{i_1i_2\cdots i_m})$ (resp. $x_i$) be the hypermatrix of
$\A$ (resp. $\x$), let
$$F_i(x)=
\sum\limits_{i_2,\cdots,i_m=1}^n a_{ii_2\cdots i_m}x_{i_2}\cdots
x_{i_m},$$ we use the notation ${\rm Res}_\x\left(\A \x^{m-1}\right)
$ to denote ${\rm Res}_x\left(F_1,...,F_n\right)$, and $\psi_\A(\lambda)$ to denote $\psi_A(\lambda)$.

\begin{theorem}
${\rm Res}_\x\left(\A \x^{m-1}\right)$ is an orthonormal invariant of $\A$.
Furthermore, all the coefficients of $\psi_\A(\lambda)$ are
orthonormal invariants of $\A$.
\end{theorem}
\proof We first prove that ${\rm Res}_\x\left(\A \x^{m-1}\right)$ is
an orthonormal invariant of $\A$.

If we choose another orthonormal frame $\tilde{e}_1,...,\tilde{e}_n$,
$$\tilde{e}_i = \sum C_i^j e_j$$
 and denote
$$\tilde{F}_i(\tilde{x})= \sum\limits_{i_2,\cdots,i_m=1}^n
\tilde{a}_{ii_2\cdots i_m}\tilde{x}_{i_2}\cdots \tilde{x}_{i_m},$$
by a direct calculation we get
$$\tilde{F}_i(\tilde{x}) = \sum C_i^jF_j(x).$$

So if we let $$G_i(x)=\sum C_i^jF_j(x),$$ then by Proposition \ref {res02},
$${\rm Res}_{x}\left(G_1,...,G_n\right)= det(C_i^j)^{(m-1)^{n-1}}{\rm Res}_x\left(F_1,...,F_n\right),$$
and noticing
$$\tilde{F}_i(\tilde{x})=G_i(\sum_j C_j^1 \tilde{x}_j, \cdots,\sum_j C_j^n \tilde{x}_j),
$$
by Corollary of Proposition \ref {res03}, we have
$${\rm Res}_{\tilde{x}}\left(\tilde{F}_1,...,\tilde{F}_n\right)= det(C_i^j)^{(m-1)^n}{\rm Res}_x\left(G_1,...,G_n\right).$$
Since $det(C_i^j)=\pm 1$, $(m-1)^{n-1}m$ is even, we have
$${\rm Res}_{\tilde{x}}\left(\tilde{F}_1,...,\tilde{F}_n\right)= {\rm Res}_x\left(F_1,...,F_n\right).$$

By the same method we can prove that when $m$ is even, ${\rm
Res}_\x\left(\A \x^{m-1}\right.$ $\left. - \lambda \left(\x^T \x\right)^{{m-2 \over 2}}\x\right)$ is orthonormal invariant, when $m$ is odd,
${\rm Res}_{(\x, x_0)}\left({\A \x^{m-1} - \lambda x_0^{m-2}\x \atop
\x^T \x - x_0^2}\right)$ is orthonormal invariant. The proof is completed.
\qed

\bigskip

Because of this theorem, in further discussion, we do not
distinguish a tensor and its corresponding hypermatrix strictly.
\begin {lemma} For a hypermatrix $A$, ${\rm Res}_x\left(A x^{m-1}\right)$ is an irreducible polynomial in the entries of $A$ (view them as variables).
\end{lemma}
\proof For any index $\alpha =(\alpha _1, \cdots , \alpha _n)$ with $| \alpha |=m-1$, it defines a subset $S_\alpha$ of $\{1, \cdots, n\}^{m-1}$,
$$S_\alpha =\{(i_2, \cdots , i_m)\ | {\rm among}\ i_2, \cdots , i_m, \ \alpha _j\ {\rm of\ them\ are\ equal\ to\ }j, \ {\rm for }\ j =1,\cdots, n\}.
$$
Obviously,
$$S_\alpha \cap S_\beta=\phi,\ {\rm for \ any}\ \alpha \not =\beta,
$$
$$\cup _\alpha S_\alpha =\{1, \cdots, n\}^{m-1}.
$$

For the system $Ax^{m-1}=0$, if we denote the equations by
$$F_i=\sum c_{i\alpha }x^\alpha,
$$
$i=1,\cdots, n$, then
$$c_{i\alpha}=\sum _{(i_2, \cdots, i_m)\in S_\alpha} a_{ii_2,\cdots, i_m}.$$

Assume ${\rm Res}_x\left(A x^{m-1}\right)=P(c_{i\alpha})=f(a_{i_1\cdots i_m})g(a_{i_1\cdots i_m})$. Let us fix a function
$$\phi: \{\alpha | \alpha=(\alpha _1, \cdots , \alpha _n), \ |\alpha |=m-1\}\to \{1, \cdots, n\}^{m-1},
$$
such that $\phi(\alpha)\in S_\alpha$.  We can choose entries of $A$
such that for every $c_{i\alpha}$ only the one entry nonzero is
$a_{i\phi(\alpha)}$. Then ${\rm Res}_x\left(A x^{m-1}\right)$ is an
irreducible polynomial in the variables corresponding to these
nonzero entries.  Hence,
$$P(a_{i\phi(\alpha)})=f(a_{i\phi(\alpha)})g(a_{i\phi(\alpha)}).
$$
Since $P$ is irreducible, one of $f(a_{i\phi(\alpha)})$,
$g(a_{i\phi(\alpha)})$, say $g(a_{i\phi(\alpha)})$ has to be of
degree 0 and the degree of the other equals to the degree of $P$.
Then $g(a_{i_1\cdots i_m})$ has degree 0 because $P(a_{i_1\cdots
i_m})$ has the same degree as $P(c_{i\alpha})$. Therefore ${\rm
Res}_x\left(A x^{m-1}\right)$ is irreducible.    \qed

We now give a general formula for $a_0(m, n)$.

\begin{theorem} \label{t2}  When $m$ is even, we have
\begin{equation} \label{a01}
a_0(m, n) = {\rm Res}_x\left(A x^{m-1}\right).
\end{equation}
When $m$ is odd, we have
\begin{equation} \label{a02}
a_0(m, n) = \left({\rm Res}_x\left(A x^{m-1}\right)\right)^2.
\end{equation}
\end{theorem}
\proof
 When $m$ is even, by (\ref{ch1}), we have
$$a_0(m, n) = \psi_A(0) = {\rm Res}_x\left(A x^{m-1}\right).$$
This proves (\ref{a01}).

When $m$ is odd, by (\ref{ch2}), we have
$$a_0(m, n) = \psi_A(0) = {\rm Res}_{(x, x_0)}\left({A x^{m-1} \atop x^Tx -
x_0^2}\right).$$ To prove (\ref{a02}), it suffices now to prove
\begin{equation} \label{a03}
{\rm Res}_{(x, x_0)}\left({A x^{m-1} \atop x^Tx - x_0^2}\right)
= \left({\rm Res}_x\left(A x^{m-1}\right)\right)^2.
\end{equation}

We first show that $\left({\rm Res}_x\left(A
x^{m-1}\right)\right)^2$ has the correct degree.   Denote $\left(A
x^{m-1}\right)_i$ as $F_i(x)$ for $i = 1 \cdots n$.   Since $x^Tx - x_0^2$ has degree $2$, by Proposition \ref{t01}, for any $i$,
${\rm Res}_{(x, x_0)}\left({A x^{m-1} \atop x^Tx -
x_0^2}\right)$ is a homogeneous polynomial in the coefficients of
$F_i(x)$, with degree $2(m-1)^{n-1}$.    Obviously, for any $i$,
$\left({\rm Res}_x\left(A x^{m-1}\right)\right)^2$ is a homogeneous
polynomial in the coefficients of $F_i(x)$, with degree
$2(m-1)^{n-1}$.   This shows that $\left({\rm Res}_x\left(A
x^{m-1}\right)\right)^2$ has the correct degree.

Next we prove that the system
\begin{equation} \label{ho1}
\left\{ {A x^{m-1} = 0, \atop x^Tx = x_0^2,} \right.
\end{equation}
has a nonzero solution if and only if ${\rm Res}_x\left(A
x^{m-1}\right)= 0$.   Let $(x, x_0)$ be a nonzero solution
of (\ref{ho1}).   Since $A$ is regular, $x_0 \not = 0$.  By the last
equation of (\ref{ho1}), $x \not = 0$.   Then $x$ is a nonzero
solution of
\begin{equation} \label{ho2}
A x^{m-1} = 0.
\end{equation}
Thus, ${\rm Res}_x\left(A x^{m-1}\right)= 0$. On the other hand,
suppose that ${\rm Res}_x\left(A x^{m-1}\right)$ $= 0$.   Then
(\ref{ho2}) has a nonzero solution $x$.  Let $x_0^2 = x^Tx$.
Then $(x, x_0)$ is a nonzero solution of (\ref{ho1}).  Therefore
${\rm Res}_{(x, x_0)}\left({A x^{m-1} \atop x^Tx - x_0^2}\right)
=0$ and ${\rm Res}_x\left(A x^{m-1}\right)$ $=0$ define the same
varieties. By the irreducibility of ${\rm Res}_x\left(A
x^{m-1}\right)$, we know
$${\rm Res}_{(x, x_0)}\left({A x^{m-1} \atop x^Tx - x_0^2}\right)
= \left({\rm Res}_x\left(A x^{m-1}\right)\right)^k.
$$
Then by degree, we know $k=2$. These prove (\ref{a03}).   The proof
is completed.
 \qed

This shows the importance of ${\rm Res}_x\left(A x^{m-1}\right)$.
When $m = 2$, it is the determinant of square matrix $A$ in the
classical sense.    Hence, it is a genuine extension of the
determinant of a square matrix, and deserves to be studied further.

\medskip

Before ending this section, we give a proposition on the degree of
the coefficients of $\psi_A(\lambda)$ as polynomials in the entries
of $A$.

When $m$ is even, we may see that
\begin{equation} \label{eig3}
A x^{m-1} - \lambda \left(x^Tx\right)^{{m-2 \over 2}}x = 0
\end{equation}
is a system of homogeneous polynomials in $x$. Every equation of
(\ref{eig3}) has the same degree $m-1$.   Thus, by Proposition
\ref{t01}, $\psi_A(\lambda)$ is a homogeneous polynomial in the
entries of $A$ and $\lambda$, with degree $n(m-1)^{n-1}$. Hence, in
(\ref{ch1}), $a_j(m,n)$ is a homogeneous polynomial in the entries
of $A$, with degree $n(m-1)^{n-1}-j$. When $m$ is odd, the first $m$
equations of
\begin{equation} \label{eig4}
\left\{{A x^{m-1} = \lambda x_0^{m-2}x \atop x^Tx =x_0^2}\right.
\end{equation}
have the same degree $m-1$, while the coefficients of last equation
of (\ref{eig4}) are either $1$ or $-1$.   Thus, by Proposition
\ref{t01}, $\psi_A(\lambda)$ is a homogeneous polynomial in the
entries of $A$ and $\lambda$, with degree $2n(m-1)^{n-1}$. Hence, in
(\ref{ch2}), $a_j(m,n)$ is a homogeneous polynomial in the entries
of $A$, with degree $2n(m-1)^{n-1}-2j$. We now have the following
proposition.

\begin{proposition} \label{degree}
In (\ref{ch1}), $a_j(m,n)$ is a homogeneous polynomial in the
entries of $A$, with degree $n(m-1)^{n-1}-j$.   In particular,
$a_{h(m,n)}(m,n)$ is a homogeneous polynomial in the entries of $A$,
with degree $n(m-1)^{n-1}- {(m-1)^n - 1 \over m-2}$.

In (\ref{ch2}), $a_j(m,n)$ is a homogeneous polynomial in the
entries of $A$, with degree $2n(m-1)^{n-1}-2j$.   In particular,
$a_{h(m,n)}(m,n)$ is a homogeneous polynomial in the entries of $A$,
with degree $2n(m-1)^{n-1}- {2(m-1)^n - 2 \over m-2}$.
\end{proposition}

\section{E-Characteristic Polynomials when $n=2$}
\label {s:CharPoly}

We now derive some simplified forms for $\psi_\A(\lambda)$ in the
case that $n=2$.

Let $b_j = \sum \{ a_{1i_2\cdots i_m} : $ exactly $j-1$ of $i_2,
\cdots, i_m$ are $2 \}$ for $j = 1, \cdots, m$, and $c_j = \sum \{
a_{2i_2\cdots i_m} : $ exactly $j-1$ of $i_2, \cdots, i_m$ are $2
\}$ for $j = 1, \cdots, m$.
For $j = 1, \cdots, m-1$, let $d_j = b_j - c_{j+1}$, and let $d_m = b_m$.

For an $N \times N$ matrix $M = (m_{ij})$, just like $\{m_{11}, m_{22}, \cdots, m_{NN} \}$ is its diagonal, we call $\{
m_{1, 2k+1}, m_{2, 2k+2}, \cdots, m_{N-2k, N}\}$ the $k$th even
upper sub-diagonal of $M$ for positive $k$  such that $2k+1 \le N$.

We now discuss the case that $m$ is even.

\begin{theorem} \label{simpeven}
Suppose that $m = 2k+2$ and $n =2$, where $k \ge 0$.  Let
$$\bar b_{2j+1} = b_{2j+1} - \left({k \atop j }\right)\lambda\ \ {\rm and}\  \bar c_{2j+2} = c_{2j+2} - \left({k \atop j }\right)\lambda,$$
for $j= 0, \cdots, k$.   Then for a regular $A$, $\psi_A(\lambda)$ is the determinant
of the following $(2m-2) \times (2m-2)$ matrix:
\begin{equation} \label{even}
M_1 = \left(\begin{array} {cccccccccc} \bar b_1 & b_2 & \bar b_3 &
\cdots & \bar b_{m-1} & b_m & 0 & 0 & \cdots &
0 \\ 0 & \bar b_1 & b_2 & \cdots & b_{m-2} & \bar b_{m-1} & b_m & 0 & \cdots & 0 \\
\cdot & \cdot & \cdot & \cdots & \cdot & \cdot & \cdot & \cdot & \cdots & \cdot \\
\vdots & \vdots & \vdots & \ddots & \vdots & \vdots & \vdots & \vdots & \ddots & \vdots \\
0 & 0 & 0 & \cdots & \bar b_1 & b_2 & \bar b_3 & b_4 & \cdots & b_m \\
0 & 0 & 0 & \cdots & c_1 & \bar c_2 & c_3 & \bar c_4 & \cdots & \bar c_m \\
-c_1 & d_1 & d_2 & \cdots & d_{m-2} & d_{m-1} & b_m  & 0 & \cdots & 0 \\
0 & -c_1 & d_1 & \cdots & d_{m-3} & d_{m-2} & d_{m-1} & b_m & \cdots & 0 \\
\vdots & \vdots & \vdots & \ddots & \vdots & \vdots & \vdots & \vdots & \ddots & \vdots \\
0 & 0 & 0 & \cdots  & d_1 & d_2 & d_3 & d_4 & \cdots & b_m \\
\end{array}\right).
\end{equation}
Here, the first $m-1$ entries of the diagonal and the first $m-2
\over 2$ even upper sub-diagonals of $M_1$ are $\bar b_1, \bar b_3,
\cdots, \bar b_{m-1}$, while the $m$th entries of the diagonal and
the first $m-2 \over 2$ even upper sub-diagonals of $M_1$ are $\bar
c_2, \bar c_4, \cdots, \bar c_m$.  They are linear factors of
$\lambda$.
\end{theorem}
\proof Consider
\begin{equation} \label{F}
F(x, \lambda) \equiv \left({\left(A x^{m-1}\right)_1 - \lambda
\left(x^Tx\right)^{{m-2 \over 2}}x_1 \atop x_2 \left(A
x^{m-1}\right)_1 - x_1\left(A x^{m-1}\right)_2 }\right) = 0.
\end{equation}
Since $A$ is regular, any nonzero solution of
 (\ref{eig3}) is a nonzero solution of (\ref{F}).  This implies that
$\psi_A(\lambda)$ is a factor of Res$_x F(x, \lambda)$.  On the
other hand, the only possible additional nonzero solution of
(\ref{F}) satisfies $x_1 = 0$ and $\left(A x^{m-1}\right)_1 = 0$.
Since $\left(A x^{m-1}\right)_1 = \sum_{i=1}^m
b_ix_1^{m-i}x_2^{i-1}$, $x_1 = 0$, $x_2 \not = 0$ and $\left(A
x^{m-1}\right)_1 = 0$ imply that $b_m = a_{12\cdots 2} = 0$.
Therefore, we conclude that
\begin{equation} \label{resF}
{\rm Res}_x F(x, \lambda) = b_m\psi_A (\lambda).
\end{equation}
Now, we have
$$F(x, \lambda) = \left({\sum_{j=0}^k \left[\bar
b_{2j+1}x_1^{m-2j-1}x_2^{2j} + b_{2j+2}x_1^{m-2j-2}x_2^{2j+1}\right]
\atop -c_1x_1^m + \sum_{i=1}^{m-1}d_ix_1^{m-i}x_2^i + d_mx_2^m
}\right)$$ By the Sylvester formula \cite{clo, gkz}, Res$_x F(x,
\lambda)$ is the determinant of the following $(2m-1) \times (2m-1)$
matrix:
$$
\left(\begin{array} {ccccccccccc} \bar b_1 & b_2 & \bar b_3 & \cdots
& \bar b_{m-1} & b_m & 0 & 0 & \cdots &
0 & 0 \\ 0 & \bar b_1 & b_2 & \cdots & b_{m-2} & \bar b_{m-1} & b_m & 0 & \cdots & 0 & 0 \\
\cdot & \cdot & \cdot & \cdots & \cdot & \cdot & \cdot & \cdot & \cdots & \cdot & \cdot \\
\vdots & \vdots & \vdots & \ddots & \vdots & \vdots & \vdots & \vdots & \ddots & \vdots & \vdots \\
0 & 0 & 0 & \cdots & \bar b_1 & b_2 & \bar b_3 & b_4 & \cdots & b_m & 0 \\
0 & 0 & 0 & \cdots & 0 & \bar b_1 & b_2 & \bar b_3 & \cdots & \bar b_{m-1} & b_m \\
-c_1 & d_1 & d_2 & \cdots & d_{m-2} & d_{m-1} & b_m  & 0 & \cdots & 0 & 0  \\
0 & -c_1 & d_1 & \cdots & d_{m-3} & d_{m-2} & d_{m-1} & b_m & \cdots & 0 & 0 \\
\vdots & \vdots & \vdots & \ddots & \vdots & \vdots & \vdots & \vdots & \ddots & \vdots & \vdots \\
0 & 0 & 0 & \cdots  & d_1 & d_2 & d_3 & d_4 & \cdots & b_m & 0 \\
0 & 0 & 0 & \cdots  & -c_1 & d_1 & d_2 & d_3 & \cdots & d_{m-1} & b_m \\
\end{array}\right).
$$
As elementary row and column operations preserve the value of the
determinant, we may subtract the last row from the $m$th row in the
above matrix. Then all the elements of the last column of the matrix
are zero except the bottom element is $b_m$.   Then we may delete
the last column and the last row, and extract $b_m$ from the
determinant. Then we see that Res$_x F(x, \lambda) =$
$b_m$det$(M_1)$. Comparing with (\ref{resF}), we have $\psi_A
(\lambda) =$ det$(M_1)$.  The proof is completed.
 \qed

The merit of $M_1$ is that only the first $m$ entries of the
diagonal and the first $m-2 \over 2$ even upper sub-diagonals of
$M_1$ contain linear factors of $\lambda$, i.e., those $\bar b_i$
and $\bar c_i$. By some elementary column operations, we may
eliminate $\lambda$ in the terms other than $\bar b_1$ and $\bar
c_2$.

When $m = 4$ and $n = 2$, we have $h(4, 2) = 4$,
$$a_{h(4, 2)=4}(4, 2) = (b_1-c_2-b_3+c_4)^2+(c_1+b_2-c_3-b_4)^2$$
$$= \left(a_{1111} + a_{2222} - a_{1221} -a_{1212} -
a_{1122} - a_{2211} - a_{2121} - a_{2112}\right)^2$$ $$+ \left(
a_{1211} + a_{1121} + a_{1112} + a_{2111} - a_{1222} - a_{2221} -
a_{2212} - a_{2122}\right)^2.$$

\medskip

We now discuss the case that $m$ is odd.  Let $e_1 = b_1c_1, e_2 =
b_1c_2 + b_2c_1, e_3 = b_1c_3 + b_2c_2 + b_3c_1, \cdots, e_m =
b_1c_m + b_2c_{m-1} + \cdots + b_mc_1, e_{m+1} = b_2c_m + \cdots +
b_mc_2, \cdots, e_{2m-1}=b_mc_m$.

\begin{theorem} \label{simpodd}
Suppose that $m = k+2$ and $n =2$, where $k \ge 1$ is odd.  Let
$$\bar e_{2j} = e_{2j} - \left({k \atop j }\right)\lambda^2,$$
for $j= 1, \cdots, m-1$.   Let $\bar f_i = \bar e_i + b_1d_{i+1}$
for $i = 2, 4, \cdots, m-1$; $f_i = e_i + b_1d_{i+1}$ for $i= 3, 5,
\cdots, m-2$; $f_i = e_i - c_md_{i-m+1}$ for $i = m, m+2, \cdots,
2m-3$; $\bar f_i = \bar e_i - c_md_{i-m+1}$ for $i = m+1, m+3,
\cdots, 2m-2$.   Then $\psi_A(\lambda)$ is the determinant of the
following $(3m-4) \times (3m-4)$ matrix:
\begin{equation} \label{odd}
M_2 = \left(\begin{array} {cccccccccc} \bar f_2 & f_3 & \cdots &
\bar f_{m-1} & e_m & \cdots & \bar e_{2m-2} & e_{2m-1} & \cdots
& 0 \\ e_1 & \bar e_2 & \cdots & e_{m-2} & \bar e_{m-1} & \cdots & e_{2m-3} & \bar e_{2m-2} & \cdots & 0 \\
\vdots & \vdots & \ddots & \vdots & \vdots & \ddots & \vdots & \vdots & \ddots & \vdots \\
0 & 0 & \cdots & \bar e_2 & e_3 & \cdots  & \bar f_{m-1} & f_m & \cdots & \bar f_{2m-2} \\
-c_1 & d_1 & \cdots & d_{m-1} & b_m & \cdots & 0 & 0 & \cdots & 0 \\
0 & -c_1 & \cdots & d_{m-2} & d_{m-1} & \cdots & 0 & 0 & \cdots & 0 \\
\vdots & \vdots & \ddots & \vdots & \vdots & \ddots & \vdots & \vdots & \ddots & \vdots \\
\cdot & \cdot & \cdots & \cdot & \cdot & \cdots & \cdot & \cdot & \cdots & \cdot \\
0 & 0 & \cdots  & 0 & 0 & \cdots &  d_2 & d_3 & \cdots & 0 \\
0 & 0 & \cdots  & 0 & 0 & \cdots & d_1 & d_2 & \cdots & b_m
\end{array}\right).
\end{equation}
Here, the first $m-2$ entries of the first row and the last $m-2$
entries of the $m$th row are somewhat different from the other
entries of the first $m$ rows, with $e_i$ or $\bar e_i$ being
replaced by $f_i$ or $\bar f_i$.
\end{theorem}
\proof Consider
\begin{equation} \label{G}
G(x, \lambda) \equiv \left({\left(A x^{m-1}\right)_1\left(A
x^{m-1}\right)_2 - \lambda^2 \left(x^Tx\right)^{{2m-4 \over
2}}x_1x_2 \atop x_2 \left(A x^{m-1}\right)_1 - x_1\left(A
x^{m-1}\right)_2 }\right) = 0.
\end{equation}
 Similar to the proof of Theorem \ref{simpeven}, we may conclude that
\begin{equation} \label{resG}
{\rm Res}_x G(x, \lambda) = b_mc_1\psi_A (\lambda).
\end{equation}
Now, we have
$$G(x, \lambda) = \left({\sum_{j=0}^{m-1}
e_{2j+1}x_1^{2m-2j-2}x_2^{2j} + \sum_{j=1}^{m-1} \bar
e_{2j}x_1^{2m-2j-1}x_2^{2j-1} \atop -c_1x_1^m +
\sum_{i=1}^{m-1}d_ix_1^{m-i}x_2^i + d_mx_2^m }\right)$$ By the
Sylvester formula \cite{clo, gkz}, Res$_x G(x, \lambda)$ is the
determinant of the following $(3m-2) \times (3m-2)$ matrix:
$$
\left(\begin{array} {cccccccccccc} e_1 & \bar e_2 & e_3 & \cdots &
\bar e_{m-1} & e_m & \cdots & \bar e_{2m-2} & e_{2m-1} & \cdots &
0 & 0 \\ 0 & e_1 & \bar e_2 & \cdots & e_{m-2} & \bar e_{m-1} & \cdots & e_{2m-3} & \bar e_{2m-2} & \cdots & 0 & 0 \\
\vdots & \vdots & \vdots & \ddots & \vdots & \vdots & \ddots & \vdots & \vdots & \ddots & \vdots & \vdots \\
0 & 0 & 0 & \cdots & \bar e_2 & e_3 & \cdots  & \bar e_{m-1} & e_m & \cdots & \bar e_{2m-2} & e_{2m-1} \\
-c_1 & d_1 & d_2 & \cdots & b_m & 0 & \cdots  & 0  & 0 & \cdots & 0 & 0  \\
0 & -c_1 & d_1 & \cdots & d_{m-1} & b_m & \cdots & 0 & 0 & \cdots & 0 & 0 \\
0 & 0 & -c_1 & \cdots & d_{m-2} & d_{m-1} & \cdots & 0 & 0 & \cdots & 0 & 0 \\
\vdots & \vdots & \vdots & \ddots & \vdots & \vdots & \ddots & \vdots & \vdots & \ddots & \vdots & \vdots \\
\cdot & \cdot & \cdot & \cdots & \cdot & \cdot & \cdots & \cdot & \cdot & \cdots & \cdot & \cdot \\
0 & 0 & 0 & \cdots  & 0 & 0 & \cdots &  d_2 & d_3 & \cdots & 0 & 0 \\
0 & 0 & 0 & \cdots  & 0 & 0 & \cdots & d_1 & d_2 & \cdots & b_m & 0 \\
0 & 0 & 0 & \cdots  & 0 & 0 & \cdots & -c_1 & d_1 & \cdots & d_{m-1} & b_m \\
\end{array}\right).
$$
Note that $e_1 = b_1c_1$ and $e_{2m-1} = b_mc_m$.   Thus, as
elementary row and column operations preserve the value of the
determinant, we may eliminate $e_1$ in the first column with a
multiple of the (m+1)th row, and eliminate $e_{2m-1}$ in the last
column with a multiple of the last row. Then we may delete the first
column, the last column, the $(m+1)$th row and the last row, and
extract $b_mc_1$ from the determinant. We see that Res$_x G(x,
\lambda) =$ $b_mc_1$det$(M_2)$. Comparing with (\ref{resG}), we have
$\psi_A (\lambda) =$ det$(M_2)$.  The proof is completed.
 \qed

The merit of $M_2$ is that only the first $m$ entries of the
diagonal and the first $m-2$ even upper sub-diagonals of $M_2$
contain linear factors of $\lambda^2$, i.e., those $\bar e_i$ and
$\bar f_i$. By some elementary column operations, we may eliminate
$\lambda$ in the terms other than $\bar f_2$ and $\bar e_2$.   


When $m = 3$ and $n = 2$, in \cite{cs}, it was given that
$$a_{h(3, 2) = 3}(3, 2) = - \left(-a_{111} + a_{122} + a_{212} +
a_{221}\right)^2 - \left(a_{112} + a_{121} + a_{211} -
a_{222}\right)^2.$$ If we use $b_i$ and $c_i$ to write it, then we
have
$$a_{h(3, 2) = 3}(3, 2) = - (b_1-c_2-b_3)^2-(c_1+b_2-c_3)^2.$$
We see that in the form using $b_i$ and $c_i$, $a_{h(3,2)}(3, 2)$
and $a_{h(4, 2)}(4, 2)$ are very similar.  Let $P_m$ and $Q_m$ be
the sum of the first $m$ terms of the following two series
respectively:
$$b_1-c_2-b_3+c_4+b_5-c_6-b_7+c_8+b_9-c_{10}- \cdots$$
and
$$c_1+b_2-c_3-b_4+c_5+b_6-c_7-b_8+c_9+b_{10}- \cdots,$$
where the signs of these terms are cyclically changed with the cycle
four.   Then, we may write
$$a_{h(3, 2)}(3, 2) = - (P_3^2+Q_3^2)$$
and
$$a_{h(4, 2)}(4, 2)= P_4^2+Q_4^2.$$
Using Theorem \ref{simpeven}, we find via computation that
$$a_{h(6, 2)}(6, 2)= (P_6^2+Q_6^2)^2.$$
It seems that there are formulas for $a_{h(m, 2)}(m, 2)$ via
$P_m^2+Q_m^2$ for all $m \ge 2$.   Later we will show that this is
true.

\section {A Bound on the Codimension of the Variety of Tensors with Infinitely Many Eigenpairs}

For positive integers $d_1, d_2, \cdots, d_n$, let $\PP$ be the complex projective space corresponding to homogeneous coordinate ring
$$\C [u_{i\alpha}\ |\ |\alpha |=\alpha _0+\cdots +\alpha _n=d_i, i=1,\cdots n].$$

Now for a system of homogeneous polynomial equations: $F_1=F_2= \cdots= F_n=0$ where $F_1, F_2, \cdots, F_n$ are homogeneous polynomials of positive degrees $d_1, d_2, \cdots, d_n$ in variables $x_0, x_1, \cdots, x_n$,
$$F_i=\sum _{|\alpha |=d_i} c_{i\alpha }x^\alpha.
$$
Then modulo obvious scalar multiplication of systems, there is one-to-one correspondence between such systems and points in $\PP$.

For such a system, the solution set in $\C \PP ^n$ has dimension $\ge n-n=0$. For systems corresponding to generic points in $\PP$, the solution sets have dimension 0, i.e. they have only finitely many solutions in $\C \PP ^n$.

Let $X\subset \PP$ be the set of points whose corresponding systems have infinitely many solutions in $\C \PP ^n$, i.e., $$X=\{([c_{i,\alpha }])\ | \ F_1=\cdots =F_n=0 \ {\rm has \ infinitely \ many \ solutions\  in} \ \C \PP ^n \}.$$ Obviously $codim \bar X\ge 1$ in $\PP$, where $\bar X$ is the Zariski closure of $X$.

Let us analysis the subset $X$ in more detail. We will use a kind of hidden variable argument similar to \cite {clo}, and we need some careful treatment to apply results on resultants, because the systems should have $n$ homogeneous equations in $n$ variables or $n$ inhomogeneous equations in $n-1$ variables.

Fix a general system $F_1=F_2= \cdots= F_n=0$,
$$F_i=\sum _{|\alpha |=d_i} u_{i\alpha }x^\alpha,
$$
where $\alpha =(\alpha _0, \alpha _1, \cdots, \alpha _n)$.

For $0\le s \not =t\le n$, we introduce a set of variables $y_j, j=0,\cdots \hat s, \cdots n$ (more precise notation would be $y^{s,t}_j$). Then in $F_i$, $i=1,\cdots, n$, we replace $x_j$ by $y_j$ for $j\not =s$, and replace $x_s$ by $\frac {x_s}{x_t}y_t$. Thus write $F_i$ as a homogeneous polynomial of degree $d_i$ in $y_j, \ j\not =s$ with coefficients in $\C [\frac {x_s}{x_t}]$:
$$\sum_{|\beta |=d_i} f_{i,\beta }^{s,t}(\frac {x_s}{x_t})y^\beta,
$$
$\beta =(\beta _0, \cdots, \hat \beta _s, \cdots, \beta_n)$, and we denote it by $F^{s,t}_i$. The set of multi-indices $\alpha $ for $x$-variables  with $\alpha _s=0$ is one-to-one corresponding to multi-indices $\beta $ for $y$-variables. Obviously
\begin {equation}
\label {eqn:st}
f^{s,t}_{i,\beta}(\frac {x_s}{x_t})=\sum _{\alpha _j=\beta _j, j\not =s,t,\alpha _s+\alpha_t=\beta _t}u_{i,\alpha}(\frac {x_s}{x_t})^{\alpha _s}.
\end{equation}

Let us make this treatment clear by an example. For the homogeneous polynomial $F=x_0^2x_1+x_0x_1^2+x_1x_2^2+x_1^2x_2$, if we take $s=0$, $t=1$, we obtain the homogeneous polynomial $[(\frac {x_0}{x_1})^2+\frac {x_0}{x_1}]y_1^3+y_1y_2^2+y_1^2y_2$ after the replacement.

The system $F_1^{s,t}=\cdots =F_n^{s,t}=0$ is a system with $n$ homogeneous equations in $n$-variables, so the resultant $Res ^{s,t}=Res (F_1^{s,t}, \cdots , F_n^{s,t})$ is a polynomial in variables $f^{s,t}_{i,\beta}$, thus a polynomial in $\Z [u_{i,\alpha}, i=1,\cdots n, \frac {x_s}{x_t}]$.

\begin {lemma} If the system $F_1=F_2=\cdots =F_n=0$ has a solution $(a_0, a_1, \cdots , a_n)$ with $a_t\not =0$, then for any $s\not= t$, $\frac {a_s}{a_t}$ is a solution of $Res (F_1^{s,t}, \cdots, F_n^{s,t})=0$.
\end{lemma}

\proof The system $\sum_{|\alpha |=d_i} f_{i,\alpha }^{s,t}(\frac
{a_s}{a_t})y^\alpha=0, i=1, \cdots, n$ has a non-trivial solution
$(\frac {a_i}{a_t})_{i\not= s}$, so $\frac {a_s}{a_t}$ is a solution
of $Res (F_1^{s,t}, \cdots, F_n^{s,t})=0$. \qed

\vskip 3mm

\begin {lemma} For any $(c_{i,\alpha })\in X$, at least one of the polynomial $Res ^{s,t}(c_{i,\alpha}, \frac {x_s}{x_t})$ vanished identically.
\end{lemma}

\proof Otherwise, a solution $(a_1,\cdots, a_n)$ of the syetem
$F_1=F_2=\cdots =F_n=0$ must satisfy $Res^{s,t}(c_{i,\alpha }, \frac
{a_s}{a_t})=0$ for any $t$ with $a_t\not =0$ and $s\not =t$, such solutions can have at most finitely many. \qed

\begin {theorem}
\label {th:codim}
The Zariski closure $\bar X$ of $X $ has codimension at least 2 in $\PP $ if $n\ge 2$.
\end{theorem}
\proof For $0\le s \not =t\le n$, let
$$X^{s,t}=\left\{[c_{i,\alpha }]\ | \ Res ^{s,t}\left(c_{i, \alpha }, \frac {x_s}{x_t}\right)\equiv 0\right\},$$
obviously $\bar X\subset \cup X^{s,t}$, so if we can prove that every $X^{s,t}$ has codimension at least 2, then so is $\bar X$.

As a polynomial in variable $\frac {x_s}{x_t}$ with coefficients in $\Z [u_{i,\alpha}, i=1,\cdots n]$, the constant term of $Res ^{s,t}(u_{i, \alpha }, \frac {x_s}{x_t})$ is the resultant for the reduced system of our general system by letting $x_s=0$, which is a general system in $n$ variables $y_j,j\not =s$. We denote the constant term by
$Res(u_{i, \alpha}, \alpha _s=0)$, by formula (\ref {eqn:st}), it is an irreducible polynomial in variables $u_{i, \alpha}$ (with $\alpha _s=0$),

Now let us consider the coefficient of linear term in $Res ^{s,t}(u_{i, \alpha }, \frac {x_s}{x_t})$. Notice by formula (\ref {eqn:st}), if we take the degree $\frac {x_s}{x_t}$ to be 0, then $Res ^{s,t}(u_{i, \alpha }, \frac {x_s}{x_t})$ is a homogeneous polynomial.

For the coefficient of linear term, again by formula (\ref {eqn:st}), every monomial is a monomial in variables $u_{i, \alpha}$ with $\alpha _s=0$ times a variable $u_{i,\alpha}$ where $\alpha _s=1$, so it has lower degree in
variables $u_{i, \alpha}$ with $\alpha _s=0$, and it has the same degree as the constant term, so it is coprime with constant term (constant term is irreducible).

Inside $\PP$ ($dim (\PP)\ge 2$ if $n\ge 2$), the subvariety defined by the constant term and the coefficient of linear term cannot be the zero set of a single polynomial, so must have codimension $\ge 2$. The subvariety $X^{s,t}$ is contained in this variety, so $X^{s,t}$ has codimension at least 2. \qed

Up to a common scalar multiplication, the set of hypermatrices is one-to-one corresponding to points in $\C \PP ^{n^m-1}$.

\begin {theorem} In the space $\C \PP ^{n^m-1}$, the set of hypermatrices with infinitely many eigenpairs has codimension at least 2 ($n\ge 2$).
\end {theorem}

\proof Let $A$ be an $m$th order hypermatrix.   Then by a similar argument as in Theorem \ref
{th:codim}, the subset of hypermatrices $A$ for which $Ax^{m-1}=0$ has
infinitely many solutions is contained in a subvariety (we denote it by $X_0$) which is defined by two coprime polynomials, and $X_0$ has codimension $\ge 2$,
i.e., the closure of the subset for which $0$ is an eigenvalue in infinitely many
eigenpairs only is a codimension $\ge 2$ subvariety.

Now for the system
$$Ax^{m-1}=\lambda x,$$
it has infinitely many solutions with $\lambda \not =0$ if and only if
$$Ax^{m-1}=x_0^{m-2} x$$
has infinitely many solutions with $x_0 \not =0$.

For $s=1,\cdots n$, we introduce variable $y=0, \cdots, \hat y_s, \cdots, y_n$, and rewrite the
system $Ax^{m-1}=x_0^{m-2} x$: replace $x_j$ by $y_j$ for $j\not =s$, and replace $x_s$ by $\frac {x_s}{x_0}y_0$, thus we have a new system
$F^s_1=F^s_2=\cdots=F^s_n=0$, this system is
obtained from $Ax^{m-1}=0$ by a replacement $x_j$ by $y_j$ for $j\not =s$ and $x_s\mapsto \frac {x_s}{x_0}y_0$ and a
translation of coefficients: with an extra term $-y_0^{m-2}y_i$ for $i\not =s$ and $-\frac {x_s}{x_0}y_0^{m-1}$ for $i=s$.

The resultant $Res ^{s}\left(a_{i_1\cdots i_n}, \frac {x_s}{x_t}\right)$ of the system $F^s_1=F^s_2=\cdots=F^s_n=0$ is a polynomial in $\Z [a_{i_1i_2\cdots i_n}, \frac {x_s}{x_0}]$, let
$$X^s=\left\{[A]\ | \ Res ^{s}\left(a_{i_1\cdots i_n}, \frac {x_s}{x_t}\right)\equiv 0\right\}\subset \C \PP ^{n^m-1} .$$
Then by similar argument in the proof of theorem \ref {th:codim}, the
set of $A$ for which $Ax^{m-1}=x_0^{m-2} x$ has infinitely many
solutions with $x_0 \not =0$ is a subset of $\cup _{s=1}^nX^{s}$, and $X^s$ has codimension at least 2.

Thus, the set of hypermatrices with infinitely many eigenpairs is a
subset of the subvariety $X_0\cup _{s=1}^n X^{s}$.  Hence its closure has codimension $\ge 2$. \qed

\section {A Bound on the Codimension of the Variety of Irregular Tensors}

Let $$X=\{[A]|A \ {\rm is \ irregular}
\}\subset \C \PP ^{n^m-1},$$

\begin {theorem} The Zariski closure $\bar X$ of $X $ has codimension at least 2.
\end{theorem}

\proof For a tensor ${ A}$, let
$$F_i=\sum a_{ii_2\cdots i_m}x_{i_2}\cdots x_{i_m}.
$$
Then ${ A}$ is irregular if and only if
$$F_1=F_2=\cdots =F_n=x^Tx=0
$$
has a solution in $\C \PP ^{n-1}$, this is a system with $n+1$ homogeneous equations in $n$ variables.

Now assume ${A}$ is irregular, then for any $i=1,\cdots, n$, the system
$$F_1=F_2=\cdots=\hat{F_i}=\cdots=F_n=x^Tx=0
$$
has a nontrivial solution, therefore the resultant $\Delta _i$ of this system is 0.

For a hypermatrix $A$ with only nontrivial entries $a_{jj\cdots jj}=1$, the system
$$F_1=F_2=\cdots=\hat{F_i}=\cdots=F_n=x^Tx=0
$$
is
$$\left\{\begin {array}{l}x_j^{m-1}=0, j=1,\cdots , \hat i, \cdots, n\\
x^Tx=0\end{array}\right.,
$$
it has no solution in $\C \PP ^{n-1}$, so $\Delta _i\not \equiv 0$.

Notice for the resultant $\Delta _i$, the variables ${a}_{ii_2\cdots
i_n}$ are missing, so the g.c.d of $\Delta _1, \cdots , \Delta _n$
is 1. i.e, any component of the subvariety of $\C \PP ^{n^m-1}$
defined by $$\Delta _1= \cdots = \Delta _n=0$$ cannot be defined by
a single polynomial. Hence any component of the subvariety $\Delta
_1= \cdots = \Delta _n=0$ has codimension $\ge 2$.

Now we have the conclusion because
$$\bar X\subset \{\Delta _1=\cdots =\Delta _n=0\}.$$
\qed



\section {The Leading Coefficient}
\label{s:LeadingCoeff}

We now study the properties of the leading coefficient
$a_{h(m, n)}(m, n)$ by exploring the difference of the definitions
of E-eigenvalues and eigenpair equivalence classes.   We see that
beside those eigenpair equivalence classes which correspond to
E-eigenvalues, all the other eigenpair equivalence classes
correspond to nonzero solutions of
\begin{equation} \label{irr}
\left\{{A x^{m-1} = \lambda x, \atop x^Tx = 0.}\right.
\end{equation}
We call (\ref{irr}) the {\bf deficit system} of $A$.

Now we investigate the leading coefficient for the characteristic polynomials for tensors when $n=2$.

\begin {lemma}
\label {lem:regular1}
For a regular tensor, the leading coefficient of its characteristic polynomial is 0, if and only if the system
$$\left \{\begin {array}{l}{A}x^{m-1}=\lambda x\\
x^Tx=0
\end{array}
\right.
$$ has a nontrivial solution.
\end{lemma}

\proof Consider the inhomogeneous system $$\left \{\begin {array}{l}{A}x^{m-1}=\lambda x,\\
x^Tx=1.
\end{array}
\right.
$$
This system has no solution at $\infty$ because of regularity.

If we consider $\lambda $ as a constant and homogenize the system with respect to $x_0$, we have the following system
$$\left \{\begin {array}{c}{A}x^{m-1}-\lambda x_0^{m-2}x=0,\\x^Tx-x_0^2=0.\end{array}\right .
$$
The resultant for this homogeneous system (the E-characteristic
polynomial) has the exact information of multiplicity of eigenpairs,
i.e., it is a constant times \cite {clo}
$$\Pi_{(\lambda_i, x_i)}(\lambda-\lambda _i)^{m(\lambda_i, x_i)},
$$
where the product is over distinguish eigenpairs $\{(\lambda_i, x_i)\}$, and $m(\lambda_i, x_i)$ is its multiplicity.

Therefore, if the leading term is 0, there is some eigenpair missing, which has to be a solution of $$\left \{\begin {array}{l}{A}x^{m-1}=\lambda x,\\
x^Tx=0.
\end{array}
\right.
$$

On the other hand, if the leading term is not 0, the deficit system cannot have a nontrivial solution, otherwise  ${A}$ will have more eigenpairs than expected.
\qed

\begin {lemma}
\label {lem:regular2}
For a regular tensor, the system
$$\left \{\begin {array}{l}{A}x^{m-1}=\lambda x,\\
x^Tx=0
\end{array}
\right.
$$ has a nontrivial solution if and only if $$P_m^2+Q_m^2=0,$$ where  $P_m$ and $Q_m$ are
the sum of the first $m$ terms of the following two series
respectively:
$$b_1-c_2-b_3+c_4+b_5-c_6-b_7+c_8+b_9-c_{10}- \cdots$$
and
$$c_1+b_2-c_3-b_4+c_5+b_6-c_7-b_8+c_9+b_{10}- \cdots,$$
where the signs of these terms are cyclically changed with the cycle
four, and $b_j = \sum \{ a_{1i_2\cdots i_m} : $ exactly $j-1$ of $i_2,
\cdots, i_m$ are $2 \}$ for $j = 1, \cdots, m$, and $c_j = \sum \{
a_{2i_2\cdots i_m} : $ exactly $j-1$ of $i_2, \cdots, i_m$ are $2
\}$ for $j = 1, \cdots, m$.
\end{lemma}
\proof Since $n=2$, if the deficit system (\ref{irr}) has a nontrivial
solution, then all the nontrivial solutions of the deficit system
(\ref{irr}) are nonzero multiple of $(1, \sqrt{-1})$ or $(1,
-\sqrt{-1})$. Therefore, there are at most two eigenpair equivalence
classes corresponding to nonzero solutions of the deficit system
(\ref{irr}), when $n=2$. Note that when $n=2$, $A x^{m-1} = \lambda
x$ can be written as
\begin{equation} \label{irr1}
\left\{{\sum_{i=1}^m b_i x_1^{m-i}x_2^{i-1} = \lambda x_1, \atop
\sum_{i=1}^m c_i x_1^{m-i}x_2^{i-1} = \lambda x_2.}\right.
\end{equation}
Substituting $(1, \sqrt{-1})$ to (\ref{irr1}), we have
$$\left\{\begin{array} {lcl} \sum_{i=1}^m b_i \left(\sqrt{-1}\right)^{i-1} & = & \lambda,
\\
\sum_{i=1}^m c_i \left(\sqrt{-1}\right)^{i-1} & = & \lambda
\sqrt{-1}.
\end{array}\right.$$
Eliminating $\lambda$ ($\lambda\not =0$ by regularity) , we have
$$\sum_{i=1}^m b_i \left(\sqrt{-1}\right)^{i-1} + \sum_{i=1}^m c_i
\left(\sqrt{-1}\right)^i = 0.$$ By the definitions of $P_m$ and
$Q_m$, we have $P_m = -Q_m\sqrt{-1}$.   Similarly, substituting $(1,
\sqrt{-1})$ to (\ref{irr1}) and eliminating $\lambda$, we have $P_m
= Q_m\sqrt{-1}$.   These imply $P_m^2+Q_m^2 = 0$.

On the other hand, if $P_m^2+Q_m^2 = 0$, by a similar but reverse
argument, we see that $(1, \sqrt{-1})$ or $(1, -\sqrt{-1})$ is a
nontrivial solution of the deficit system (\ref{irr}).   This completes
the proof.  \qed

\begin {theorem}
\label{th:PQ} If ${A}$ has only finitely many equivalence
classes of eigenpairs, then the leading coefficient
$a_{h(m,2)}(m,2)=0$ if and only if
$$P_m^2+Q_m^2=0.$$
\end{theorem}
\proof First, if ${A}$ has only finitely many equivalence classes of
eigenpairs, then its characteristic polynomial cannot be identically
0.

Let us consider the subvarieties of $\C \PP ^{2^m-1}$:
$X=\{a_{h(m,2)}(m,2)=0\}$, and $Y=\{P_m^2+Q_m^2=0\}$, they both have
codimension 1. Let $Z_1=\{[A]|{\rm \ A\ has\ infinitely\ many\
eigenpairs}\}$, and $Z_2=\{[A]|\ A{\rm \ is\ irregular}\}$. Then
both $Z_1$ and $Z_2$ have codimension $\ge 2$. So $X-\bar Z_1-\bar
Z_2$ and $Y-\bar Z_1-\bar Z_2$ are not empty.

By Lemmas \ref {lem:regular1} and \ref {lem:regular2}, $X-\bar Z_1-\bar Z_2\subset Y$ and $Y-\bar Z_1-\bar Z_2\subset X$, taking closure, we have $X\subset Y$ and $Y \subset X$. This implies $X=Y$ and the conclusion follows. \qed

\begin {theorem} If ${ A}$ has only finitely many equivalence classes of eigenpairs, then the leading coefficient $a_{h(m,2)}(m,2)$ is
$$(P_m^2+Q_m^2)^{\frac {m-2}2}$$
for even $m$, and
$$-(P_m^2+Q_m^2)^{{m-2}}$$
for odd $m$.
\end{theorem}

\proof Denote that
$$W = \{ [B] | P_m^2 + Q_m^2 = 0 \}.$$
$W$
can be written as a union of two irreducible varieties
$$W = W_1 \bigcup W_2,$$
where
$$W_1=\{ [B] | P_m+\sqrt{-1}Q_m=0\},$$
$$W_2=\{ [B] | P_m-\sqrt{-1}Q_m=0\}.$$
Since $P_m+\sqrt{-1}Q_m$ and $P_m-\sqrt{-1}Q_m$ are irreducible, by
Theorem \ref {th:PQ}
$$a_{h(m, 2)}(m, 2)=(P_m+\sqrt{-1}Q_m)^k(P_m-\sqrt{-1}Q_m)^sC,$$
where $C$ is a constant. Since  $a_{h(m, 2)}(m, 2) \in
\Z[a_{i_1\cdots i_m}]$, $k$ must be equal to $s$. Then by checking
the degree, we see that our results hold up to a scaling constant. By
checking an example such as $a_{1\cdots 1} = a_{2\cdots 2} = 1$ but
all the other entries of $A$ are zero, we have the results. \qed

{\bf Acknowledgment}  We are thankful to two referees, whose
comments greatly improved our paper.

\end{document}